\documentclass[12pt]{amsart}
\usepackage{amssymb,amsmath}

\headsep=1cm \numberwithin{equation}{section}
\newtheorem{theorem}{Theorem}[section]
\newtheorem{lemma}[theorem]{Lemma}

\newtheorem{corollary}[theorem]{Corollary}

\theoremstyle{definition}
\newtheorem{definition}[theorem]{Definition}
\theoremstyle{remark}
\newtheorem{remark}[theorem]{Remark}
\newtheorem{example}[theorem]{Example}

\newcommand{\Ass}{\operatorname{Ass}}
\newcommand{\im}{\operatorname{im}}

\newcommand{\Spec}{\operatorname{Spec}}

\newcommand{\E}{\operatorname{E}}

\newcommand{\V}{\operatorname{V}}

\newcommand{\Ext}{\operatorname{Ext}}
\newcommand{\Supp}{\operatorname{Supp}}
\newcommand{\Tor}{\operatorname{Tor}}
\newcommand{\Hom}{\operatorname{Hom}}

\newcommand{\Ann}{\operatorname{Ann}}

\newcommand{\lo}{\longrightarrow}

\newcommand{\fa}{\frak{a}}

\begin{document}
\author[Divaani-Aazar and Mafi ]{Kamran Divaani-Aazar and Amir Mafi}
\title[Associated primes of local cohomology modules]
{Associated primes of local cohomology modules}

\address{K. Divaani-Aazar, Department of Mathematics, Az-Zahra University,
Vanak, Post Code 19834, Tehran, Iran and Institute for Studies in
Theoretical Physics and Mathematics, P. O. Box 19395-5746, Tehran,
Iran.} \email{kdivaani@ipm.ir}

\address{A. Mafi,  Institute of Mathematics, University for Teacher
Education, 599 Taleghani Avenue, Tehran 15614, Iran.}



\begin{abstract} Let $\fa$ be an ideal of a commutative Noetherian
ring $R$ and $M$ a finitely generated $R$-module. Let $t$ be a
natural integer. It is shown that there is a finite subset $X$ of
$\Spec R$, such that $\Ass_R(H_{\fa}^t(M))$ is contained in $X$
union with the union of the sets
$\Ass_R(\Ext_R^j(R/\fa,H_{\fa}^i(M)))$, where $0\leq i<t$ and
$0\leq j\leq t^2+1$. As an immediate consequence, we deduce that
the first non $\fa$-cofinite local cohomology module of $M$ with
respect to $\fa$ has only finitely many associated prime ideals.
\end{abstract}

\maketitle

\section{Introduction}

Throughout this paper, $R$ is a commutative Noetherian ring with
identity. For an ideal $\fa$ of $R$ and an $R$-module $M$, the
$i$-th local cohomology module of $M$ with respect to $\fa$ is
defined as:
$$H_{\fa}^i(M)=\underset{n}{\varinjlim}\Ext_R^i(R/\fa^n,M).$$
The reader can refer to [{\bf 3}], for the basic properties of local
cohomology.

In [{\bf 6}], Hartshorne defines an $R$-module $M$ to be
$\fa$-cofinite if $\Supp_RM\subseteq \V(\fa)$ and
$\Ext_R^i(R/\fa,M)$ is finitely generated for all $i\geq 0$. He
asks when the local cohomology modules of a finitely generated
module are $\fa$-cofinite. In this regard, the best known result
is that for a finitely generated $R$-module $M$ if either $\fa$ is
principal or $R$ is local and $\dim R/\fa=1$, then the modules
$H_{\fa}^i(M)$ are $\fa$-cofinite. These results are proved in
[{\bf 8}, Theorem 1] and [{\bf 14}, Theorem 1.1], respectively.

Since for an $\fa$-cofinite module $N$, we have
$\Ass_RN=\Ass_R(\Hom_R(R/\fa,N))$, it turns out that $\Ass_RN$ is
finite. Huneke [{\bf 7}] raised the following question: If $M$ is
a finitely generated $R$-module, then the set of associated primes
of $H_{\fa}^i(M)$ is finite for all ideals $\fa$ of $R$ and all
$i\geq 0$. Singh [{\bf 12}] gives a counter-example to this
conjecture. However, it is known that this conjecture is true in
many situations. For example, Brodmann and Lashgari [{\bf 2},
Theorem 2.2] showed that, if for a finitely generated $R$-module
$M$ and an integer $t$, the local cohomology modules
$H_{\fa}^0(M),H_{\fa}^1(M),\dots ,H_{\fa}^{t-1}(M)$ are all
finitely generated, then $\Ass_R(H_{\fa}^t(M))$ is  finite. For a
survey of recent developments on finiteness properties of local
cohomology, see Lyubeznik's interesting article [{\bf 10}].

In this article, we first introduce the class of weakly Laskerian
modules. This class includes all Noetherian modules and also all
Artinian modules. Moreover, this class is large enough to contain
all Matlis reflexive modules as well as all linear compact
modules. Then as the main result of this paper, we establish the
following. Let $M$ be a weakly Laskerian module and $t\in
\mathbb{N}$ a given integer. There is a finite subset $X$ of
$\Spec R$ such that
$$\Ass_R(H_{\fa}^t(M))\subseteq
\bigl(\bigcup_{0\leq i<t,0\leq j\leq t^2+1}
\Ass_R(\Ext_R^j(R/\fa,H_{\fa}^i(M)))\bigr)\cup X.$$ Clearly this
result implies the main result of [{\bf 2}].

\section{The results}

An $R$-module $M$ is said to be Laskerian if any submodule of $M$
is an intersection of a finite number of primary submodules.
Obviously, any Noetherian module is Laskerian. Next, we present
the following definition.

\begin{definition} An $R$-module $M$ is said to be {\it weakly
Laskerian} if the set of associated primes of any quotient module
of $M$ is finite.
\end{definition}

\begin{example} i) Any Laskerian module is weakly Laskerian. In
particular, any Noetherian module is weakly Laskerian.\\
ii) It is known that the set of associated primes of  an Artinian
module is a finite set consisting of maximal ideals. Hence any
Artinian module is weakly Laskerian.\\
iii) Recall that a module $M$ is said to have finite Goldie
dimension if $M$ does not contain an infinite direct sum of
non-zero submodules, or equivalently, the injective envelope
$\E(M)$ of $M$ decomposes as a finite direct sum of indecomposable
injective submodules. Because for any $R$-module $C$, we have
$\Ass_R(C)=\Ass_R(\E(C))$, it turns out that any module with
finite Goldie dimension has only finitely many associated prime
ideals. This yields that a module all of whose
quotients have finite Goldie dimension is weakly Laskerian.\\
iv) Let $E$ be the minimal injective cogenerator of $R$ and $M$ an
$R$-module. If for an $R$-module $M$ the natural map from $M$ to
$\Hom_R(\Hom_R(M,E),E)$ is an isomorphism, then $M$ is said to be
Matlis reflexive. By [{\bf 1}, Theorem 12], an $R$-module $M$ is
Matlis reflexive if and only if $M$ has a finitely generated
submodule $S$ such that $M/S$ is Artinian and $R/\Ann_RM$ is a
complete semi-local ring. Also, as it is mentioned in [{\bf 5},
Corollary 1.2], one can deduce from the argument [{\bf 4},
Proposition 1.3], that any quotient of an $R$-module $M$ has
finite Goldie dimension if and only if $M$ has a finitely
generated submodule $S$ such that $M/S$ is Artinian. Thus, by
(iii), any Matlis reflexive module is weakly Laskerian.\\
v) An $R$-module $M$ is said to be linearly compact if each system
of congruences $$x\equiv x_i (M_i),$$ indexed by a set $I$ and
where the $M_i$ are submodules of $M$, has a solution $x$ whenever
it has a solution for every finite subsystem. It is known that the
category of linearly compact $R$-modules form a Serre subcategory
of the category of all $R$-modules. In particular, every quotient
of a linearly compact module is also linearly compact. On the
other hand a linearly compact module $M$ has finite Goldie
dimension (see e.g. [{\bf 13}, Chapter 1.3]). Thus, if $M$ is a
linearly compact module, then any quotient of $M$ has finite
Goldie dimension, and so $M$ is weakly Laskerian by (iii).
\end{example}

To prove the main result of this paper, we need to the following
two lemmas.

\begin{lemma} i) Let  $0\lo L\lo M\lo N\lo 0$, be an exact
sequence of $R$-modules. Then $M$ is weakly Laskerian if and only
if $L$ and $N$ are both weakly Laskerian. Thus any subquotient of
a weakly Laskerian module as well as any finite direct sum of
weakly Laskerian modules is weakly Laskerian.\\
ii) Let $M$ and $N$ be two $R$-modules. If $M$ is weakly Laskerian
and $N$ is finitely generated, then $\Ext^i_R(N,M)$ and
$\Tor^R_i(N,M)$ are weakly Laskerian for all $i\geq 0$.
\end{lemma}

{\bf Proof.} The proof of (i) is easy and we leave it to the
reader.\\
ii) We only prove the assertion for the $\Ext$ modules and the
proof for the $\Tor$ modules is similar. Because $R$ is a
Noetherian ring and $N$ is finitely generated, it follows that $N$
possesses a free resolution $$F_{\cdot}:  \dots \lo
F_n\overset{d_n}\lo F_{n-1}\overset{d_{n-1}}\lo \dots
\overset{d_2}\lo F_1\overset{d_1}\lo F_0\lo 0,$$ consisting of
finitely generated free modules. If $F_i=\oplus^nR$ for some
integer $n$, then $\Ext^i_R(N,M)=H^i(\Hom_R(F_{\cdot},M))$ is a
subquotient of $\oplus^nM$. Therefore, it follows from (i), that
$\Ext_R^i(N,M)$ is weakly Laskerian for all $i\geq 0$. $\Box$

By [{\bf 11}, Theorem 11.38], there is a Grothendieck spectral
sequence with
$E_2^{p,q}:=\Ext_R^p(R/\fa,H^q_{\fa}(M))\underset{p}{\Longrightarrow}
\Ext_R^{p+q}(R/\fa,M)$. Let $E_{\infty}:=\{E_{\infty}^{p,q}\}$ be
the limit term of this spectral sequence. In the sequel, we show
that, if $M$ is weakly Laskerian, then $\Ass_R(E_{\infty}^{p,q})$
is finite for all $p,q$ with $0\leq p \leq q$.

\begin{lemma} Let $\fa$ be an ideal of $R$. If $M$ is a weakly
Laskerian module, then the set of associated primes of
$E_{\infty}^{p,q}$ is finite for all $p,q$ with $0\leq p \leq q$.
\end{lemma}

{\bf Proof.} Since, by [{\bf 11}, Theorem 11.38], the Grothendieck
spectral sequence $E_2^{p,q}=\Ext_R^p(R/\fa,H^q_{\fa}(M))$
converges to $H^{p+q}:=\Ext_R^{p+q}(R/\fa,M)$, it follows that
there is a finite filtration $$0=\phi^{q+1}H^q\subseteq
\phi^qH^q\subseteq \dots \subseteq \phi^1H^q \subseteq
\phi^0H^q=H^q,$$ of $H^q$ such that $E_{\infty}^{p,q}\cong
\phi^pH^q/\phi^{p+1}H^q$ for all $p=0,1,\dots ,q$. Because $M$ is
weakly Laskerian, by Lemma 2.3 (ii), it turns out that
$\Ext_R^q(R/\fa,M)$ is also weakly Laskerian. Hence any
subquotient of $H^q$ is weakly Laskerian. In particular,
$\Ass_R(E_{\infty}^{p,q})$ is finite. $\Box$

Now, we are ready to prove the main theorem of this paper.

\begin{theorem} Let $\fa$ be an ideal of $R$ and $M$ a weakly
Laskerian $R$-module. Let $t$ be a natural integer. There is a
finite subset $X$ of $\Spec R$ such that $$\Ass_R(\Ext_R^l(R/\fa,
H_{\fa}^t(M)))\subseteq \bigl(\bigcup_{0\leq i<t, 0\leq j\leq
t^2+1} \Ass_R(\Ext_R^j(R/\fa,H_{\fa}^i(M)))\bigr)\cup X,$$ for
$l=0,1$.
\end{theorem}

{\bf Proof.} Consider the Grothendieck spectral sequence
$$E_2^{p,q}:=\Ext_R^p(R/\fa,H^q_{\fa}(M))\underset{p}{\Longrightarrow}
\Ext_R^{p+q}(R/\fa,M).$$ Set
$$X=\bigl(\bigcup_{0\leq j\leq t^2+1}
\Ass_R(\Ext_R^j(R/\fa,H^0_{\fa}(M)))\bigr)\cup
\Ass_R(E_{\infty}^{0,t})\cup \Ass_R(E_{\infty}^{1,t}).$$ Then $X$
is  a finite set, by Lemmas 2.3 and 2.4. First, we prove the claim
for $l=0$. We have to show that $$\Ass_R(E_2^{0,t})\subseteq
\bigl(\bigcup_{0\leq i<t, 0\leq j\leq t^2+1}
\Ass_R(E_2^{j,i})\bigr)\cup X.$$ From the choice of $X$, it is
clear that we may assume $M$ is $\fa$-torsion free. The exact
sequence
$$0\lo \ker d_i^{0,t}\lo E_i^{0,t}\overset{d_i^{0,t}}\lo
E_i^{i,t-i+1},$$ yields that $\Ass_R(E_i^{0,t})\subseteq
\Ass_R(\ker d_i^{0,t})\cup \Ass_R(E_i^{i,t-i+1})$ for all $i\geq
2$. But $\ker d_i^{0,t}=E_{i+1}^{0,t}$, because
$E_i^{-i,t+i-1}=0$. Hence
$$\Ass_R(E_i^{0,t})\subseteq \Ass_R(E_{i+1}^{0,t})\cup
\Ass_R(E_i^{i,t-i+1}), (1)$$ for all $i\geq 2$.

Let $n>t$ be an integer and consider the sequence
$${\tiny {E_n^{-n,t+n-1}}}\overset{{\tiny {d_n^{-n,t+n-1}}}}
\longrightarrow E_n^{0,t}
\overset{d_n^{0,t}}\lo E_n^{n,t-n+1}.$$ Since $M$ is $\fa$-torsion
free, $E_n^{n,0}=0$. Note that for each $i\geq 2$, the module
$E_i^{p,q}$ is a subquotient of $E_2^{p,q}$. Also, $E_n^{i,j}=0$
if either $i< 0$ or $j< 0$. Thus, we have $\ker
d_n^{0,t}=E_n^{0,t}$ and $\im d_n^{-n,t+n-1}=0$, and so
$$E_{n+1}^{0,t}=\ker d_n^{0,t}/\im d_n^{-n,t+n-1}\cong E_n^{0,t}.
 (2)$$

Using (2) successively for all $n>t$, we get $E_{t+1}^{0,t}\cong
E_{t+2}^{0,t}\cong \dots =E_{\infty}^{0,t}$.  Now, by iterating
(1) for all $i=2,\dots ,t$, we deduce that
$$\Ass_R(E_2^{0,t})\subseteq
(\bigcup_{i=2}^t\Ass_R(E_i^{i,t-i+1}))\cup X.$$

Next, we show that
$$\Ass_R(E_i^{i,t-i+1})\subseteq
\bigcup_{k=1}^t\Ass_R(E_2^{ki,t-ki+k}),$$ for all $i=3,\dots ,t$.
Clearly this finishes the proof for the case $l=0$. Consider the
exact sequence
$$0\lo \ker d_i^{ki,t-ki+k}\lo {\tiny
{E_i^{ki,t-ki+k}}}\overset{{\tiny {
d_i^{ki,t-ki+k}}}}\longrightarrow E_i^{(k+1)i,t-(k+1)i+k+1}.$$
Since $\ker d_i^{ki,t-ki+k}\subseteq \ker d_2^{ki,t-ki+k}\subseteq
E_2^{ki,t-ki+k}$ and $E_i^{p,q}=0$ for all $q\leq 0$, by using the
above exact sequence successively for $k=1,2,\dots ,t$, we deduce
that $\Ass_R(E_i^{i,t-i+1})\subseteq
\bigcup_{k=1}^t\Ass_R(E_2^{ki,t-ki+k})$.

Now, by repeating the above argument, we can show that
$$\Ass_R(E_2^{1,t})\subseteq \bigl(\bigcup_{0\leq i<t, 0\leq j\leq
t^2+1} \Ass_R(E_2^{j,i})\bigr)\cup X.$$ Therefore the proof is
complete. $\Box$

Now, we can obtain the following extension of [{\bf 2}, Theorem
2.2]. Note that, because $\Supp_R(H_{\fa}^i(M))\subseteq \V(\fa)$,
it follows that $H_{\fa}^i(M)$ is $\fa$-cofinite, whenever it is
finitely generated.

\begin{corollary} Let $\fa$ be an ideal of $R$ and $M$ a weakly
Laskerian module. Let $t\in \mathbb{N}_0$ be an integer such that
$H_{\fa}^i(M)$ is $\fa$-cofinite for all $i< t$. Then, the sets of
associated primes of $H_{\fa}^t(M)$ and of
$\Ext_R^1(R/\fa,H_{\fa}^t(M))$ are finite.
\end{corollary}

{\bf Proof.} If $t=0$, then the claim follows by Lemma 2.3. Now
assume that $t>0$ and let $i< t$ be an integer. Because
$H_{\fa}^i(M)$ is $\fa$-cofinite for any $j\geq 0$,
$\Ext_R^j(R/\fa,H_{\fa}^i(M))$ is finitely generated, and so
$\Ass_R(\Ext_R^j(R/\fa,H_{\fa}^i(M)))$ is finite. On the other
hand, since $\Supp_R(H_{\fa}^t(M))\subseteq \V(\fa)$, we have
$$\Ass_R(\Hom_R(R/\fa,H_{\fa}^t(M)))=\Ass_R(H_{\fa}^t(M))\cap \V(\fa)
=\Ass_R(H_{\fa}^t(M)).$$ Therefore the conclusion follows by
Theorem 2.5. $\Box$

It is clear by Theorem 2.5, that $\Ass_R(H_{\fa}^t(M))$ is finite,
whenever the sets of associated primes of the modules
$\Ext_R^j(R/\fa,H_{\fa}^i(M))$  are finite for all $0\leq i<t$ and
all $0\leq j\leq t^2+1$. Thus Lemma 2.3 (ii) yields the following.

\begin{corollary} Let $\fa$ be an ideal of $R$ and $M$ a weakly
Laskerian module. Let $t\in \mathbb{N}_0$ be an integer such that
$H_{\fa}^i(M)$ is weakly Laskerian module for all  $i< t$. Then,
the sets of associated primes of $H_{\fa}^t(M)$ and of
$\Ext_R^1(R/\fa,H_{\fa}^t(M))$ are finite.
\end{corollary}

\begin{remark} Let $M$ be a finitely generated $R$-module.
Khashyarmanesh and Salarian [{\bf 9}, Theorem B($\beta$)] have
proved that if $t$ is an integer such that $\Supp_R(H_{\fa}^i(M))$
is finite for all $i<t$, then the set of associated primes of
$H_{\fa}^t(M)$ is finite. Clearly any $R$-module with finite
support is weakly Laskerian. Hence Corollary 2.7  generalizes
[{\bf 9}, Theorem B($\beta$)].
\end{remark}


\end{document}